\documentclass[a4paper,12pt]{amsart}
\usepackage[pdfpagelabels]{hyperref}

\usepackage{amsmath}
\usepackage{amsthm}
\usepackage{amssymb}
\usepackage{graphicx}
\usepackage{cases}
\usepackage{bm}
\usepackage{euscript}
\usepackage{mfirstuc,textcase}
\newcommand{\R}{\mathbb{R}}
\newcommand{\D}{\mathbb{D}}
\newcommand{\C}{\mathbb{C}}

\renewcommand{\H}{\mathcal{H}}
\renewcommand{\S}{\mathcal{S}}

\newcommand{\N}{\mathbb{N}}
\numberwithin{equation}{section}
\newcommand{\abs}[1]{\lvert#1\rvert}

      \def \ov{\overline}         \def\beq{\begin{equation}} \def\beqq{\begin{equation*}}
\def\eeq{\end{equation}}    \def\eeqq{\end{equation*}} \def\bprof{\begin{proof}}    \def\eprof{\end{proof}}  \def\bald{\begin{aligned}}    \def\eald{\end{aligned}}
         \def\bthm{\begin{thm}}
\def\ethm{\end{thm}}   \def\p{\partial}    
\theoremstyle{cupthm}
\newtheorem{thm}{Theorem}[section]

\newtheorem{lem}[thm]{Lemma}
\theoremstyle{cupdefn}

\theoremstyle{cuprem}

\numberwithin{equation}{section}
\newtheorem{innercustomthm}{Theorem}
\newenvironment{customthm}[1]
  {\renewcommand\theinnercustomthm{#1}\innercustomthm}
  {\endinnercustomthm}
\newtheorem{problem}{Problem}[section]

\begin{document}
\title[V\MakeLowercase{riability Regions for the $n$-th Derivative}]{V\MakeLowercase{ariability Regions for the $n$-th Derivative of Bounded Analytic Functions}}
\author{Gangqiang Chen}
\address{School of Mathematics and Computer Sciences,
Nanchang University, Nanchang 330031, China}
\address{Graduate School of Information Sciences, Tohoku University, Sendai 980-8579, Japan}
\email{cgqmath@qq.com; chenmath@ncu.edu.cn}
\date{\today}



\begin{abstract}
 Let  $\mathcal{H}$ be the class of all analytic self-maps of the open unit disk $\mathbb{D}$. Denote by $H^n f(z)$ the $n$-th order hyperbolic derivative of $f\in \mathcal H$ at  $z\in \mathbb{D}$.
 For $z_0\in \mathbb{D}$ and $\gamma  = (\gamma_0, \gamma_1 , \ldots , \gamma_{n-1}) \in {\mathbb D}^{n}$, let ${\mathcal H} (\gamma) = \{f \in {\mathcal H} : f (z_0) = \gamma_0,H^1f (z_0) = \gamma_1,\ldots ,H^{n-1}f (z_0) = \gamma_{n-1} \}$. In this paper, we  determine the variability region $V(z_0, \gamma ) = \{ f^{(n)}(z_0) : f \in {\mathcal H} (\gamma) \}$, which can be called ``the generalized Schwarz-Pick Lemma of $n$-th derivative". We then apply the generalized Schwarz-Pick Lemma to establish a $n$-th order Dieudonn\'e's Lemma, which provides an explicit description of the variability region
$\{h^{(n)}(z_0): h\in \mathcal{H}, h(0)=0,h(z_0) =w_0, h'(z_0)=w_1,\ldots, h^{(n-1)}(z_0)=w_{n-1}\}$ for given $z_0$, $w_0$, $w_1,\dots,w_{n-1}$.
Moreover, we determine the form of all extremal functions.
\end{abstract}


\maketitle
\section{Introduction}
\label{intro}
We denote by $\D$ and $\overline{\D}$ the open and closed unit disks in  the complex plane $\C$, and define the disks
$\D(c, \rho):=\left\{ \zeta \in \C : |\zeta-c|< \rho \right\}$,
and
$\overline{\D}(c, \rho):=\left\{\zeta \in \C : |\zeta-c|\le \rho \right\}$ for $c\in\C$ and $\rho>0$.
Throughout this article, let $z_0\in \D$ be a given point,  $\H$ be the class of all analytic self-mappings of $\mathbb{D}$ and $\H_0=\{f\in\H:f(0)=0\}$.
 In addition, $\mathcal{S}$ will denote the set of analytic functions from $\D$ to $\overline{\D}$.

First we recall the classical Schwarz's Lemma, which
says that if $f \in \mathcal{H}_0$, then $| f(z_0)| \le |z_0|$ for any non-zero $z_0$ in $\D$ and $|f'(0)|\le 1$, and equalities hold if and only if
 $ f(z)=e^{i \theta}z$ for some $\theta \in \mathbb{R}$.
Since then, more and more authors started to consider the bounded analytic functions. More recently,  a lot of articles on regions of variability have been written (for instance, see \cite{chen_2019,chen2020,chen2021,Hoshinaga-Yanagihara2024Julia,kaptanoglu2002refine,mercer1997sharpened,yanagihara2005}).
 In 1916, Pick \cite{Pick1916} proved the well-known Schwarz-Pick Lemma, which states that
 $$|f'(z)|\le \frac{1-|f(z)|^2}{1-|z|^2},\quad f\in \mathcal{H}, \quad z\in \D,$$
 and equality holds if and only if
$f(z)=e^{i\theta }(z-a)/(1-\bar a z)$, $\theta\in \R$, $a \in \D$.
If we let $T_{a}(z):=(z+a) /(1+\bar{a} z)$, $a,z\in \D$,
then the
Schwarz-Pick Lemma can be restated in a modified form as follows.
\begin{customthm}{A (Schwarz-Pick Lemma)}\label{thm:g-first}
Let $z_0, \gamma_0\in \D$.
Suppose that $f\in \H$,
$f(z_0)=\gamma_0$.
Set
$$
 f_{\alpha}(z) = T_{\gamma_0}\big(\alpha T_{-z_0}(z)\big).$$
Then the region of values of $f'(z_0)$ is the closed disk
$$
\overline{\D}(0, \frac{1-|\gamma_0|^2}{1-|z_0|^2})
=\{f'_{\alpha}(z_0):\alpha\in\overline{\D} \},
$$
and $f(z)$ is the form of
$T_{\gamma_0}\big(T_{-z_0}(z) f^*(z)\big)$,
where $f^*\in\S$.
Further, $f'(z_0)\in \partial\D(0, \dfrac{1-|\gamma_0|^2}{1-|z_0|^2})$ if and only if
$f(z)=f_{\alpha}(z)$ for some constant $\alpha \in \partial \D$.
\end{customthm}

Based on the Schwarz-Pick Lemma, Dieudonn\'e \cite{dieudonne1931} obtained the variability region of $f'(z_0)$ for $f\in \mathcal{H}_0$, which was called
Dieudonn\'e's Lemma (see also \cite{beardon2004multi} and \cite{duren1983univalent}).
\begin{customthm}{B (Dieudonn\'e's Lemma)}
Let $z_0,w_0\in \D$ with $|w_0|=s<r=|z_0|$. Suppose that $f\in\mathcal{H}_0$ and $f(z_0) = w_0$.
Set $u_0=w_0/z_0$, $f_{\alpha}(z)=zT_{u_0}(\alpha T_{-z_0}(z))$. Then the region of values of $f'(z_0)$ is the closed disk
\begin{align*}
\overline{\mathbb{D}}
  \left( c'_1(z_0,w_0) , \rho'_1(z_0,w_0) \right)=\{f'_{\alpha}(z_0):\alpha\in \overline{\D}\},
\end{align*}
where $$c'_1(z_0,w_0)=\frac{w_0}{z_0}, \qquad \rho'_1(z_0,w_0)=\frac{|z_0|^2-|w_0|^2}{|z_0|(1-|w_0|^2)},$$
and $f(z)$ is the form of
$$z T_{u_0}\left(T_{-z_0}(z) f^*(z)\right),$$
where $f^*\in \S$.
Further,
$f'(z_0) \in \partial \D \left( c'_1(z_0,w_0), \rho'_1(z_0,w_0)\right) $ if and only if
$f(z)=f_{\alpha}(z)$ for some constant $\alpha\in \partial \D$.
\end{customthm}

It is natural to describe the variability region of  $f^{(n)}(z_0)$ for $f\in\H$ and further
establish a $n$-th order Dieudonn\'e's Lemma to determine the region of values of $h^{(n)}(z_0)$ for $h\in \H_0$ not a rotation about the origin. Indeed, assume that $f(z)=h(z)/z$, then $f$ is an analytic self-mapping of $\D$. Since $h^{(n)}(z)= zf^{(n)}(z)+nf^{(n-1)}(z)$, we just need to determine the variability region of $f^{(n)}(z_0)$ for $f\in \H$, which is related to the $n$-th order hyperbolic derivative of $f$ at $z_0$.
The notion of higher-order hyperbolic derivatives for $f\in \H$ was  introduced by P. Rivard \cite{rivard2011higher-orderHyperbolicDerivatives} (see also \cite{Baribeau2013}).
For $ z, w\in\D$, we define
\begin{equation}\label{eq:dist}
[z,\,w]:=
\begin{cases}
\dfrac{z-w}{1-\overline{w}z} &\quad \text{if}~z\bar w\ne -1; \\
\infty &\quad\text{if}~z\bar w= -1.
\end{cases}
\end{equation}
We construct an operator $\Delta_{z_0}$, which maps every function $f\in \H$ to $\Delta_{z_0} f\in \S$, by
\begin{equation}\label{eq:hdquo}
\Delta_{z_0}f(z)=
\begin{cases}
\dfrac{[f(z),\,f(z_0)]}{[z,\,z_0]} &\quad \text{for}~z\ne z_0, \\
\null & \\
\dfrac{(1-|z_0|^2)f'(z_0)}{1-|f(z_0)|^2} &\quad\text{for}~z=z_0.
\end{cases}
\end{equation}
Then we can iterate the process and construct the hyperbolic divided difference of order $j$ of the function $f$ for distinct parameters $z_0,\cdots,z_{j-1}$ as follows (cf. \cite{Baribeau-Rivard2009}):
$$
\Delta^jf(z;z_{j-1},\cdots,z_0)=(\Delta_{z_{j-1}}\circ\cdots\circ\Delta_{z_0})f(z).
$$
We can rewrite the above recursive definition as
\begin{equation}
\Delta^jf(z;z_j,\cdots,z_0)=\frac{[\Delta^{j-1}f(z;z_{j-1},\cdots,z_0),\Delta^{j-1}f(z_j;z_{j-1},\cdots,z_0)]}{[z,z_j]}.
\end{equation}
Using a limiting process, we define $\Delta_z^n f(\zeta)$ by
\begin{equation}\label{def:Delta-z-f-zeta}
\Delta_z^n f(\zeta)=\Delta^n f(\zeta;z,\cdots,z):=\lim_{z_{n-1}\to z}\cdots\lim_{z_0\to z}\Delta^n f(\zeta;z_{n-1},\cdots,z_0)\quad(\zeta\neq z).
\end{equation}
Therefore, we define the $n$-th order  hyperbolic derivative  of $f\in \H$ at  $z\in \D$ by
$$H^n f(z):=\Delta_z^n f(z)=\Delta^n f(z;z,\cdots,z):=\lim_{\zeta\to z}\Delta^n f(\zeta;z,\cdots,z).$$
The usual hyperbolic derivative coincides with the first-order hyperbolic derivative $H^{1} f$,
$$f^h(z):=\frac{(1-|z|^2)f'(z)}{1-|f(z)|^2}=H^{1} f(z).$$

\section{Schwarz-Pick Lemma of the $n$-th derivative}
We define a sequence of functions $F_n$ of the $n$ complex variables
$\gamma_1, \gamma_2, \dots, \gamma_n~(n=1,2,3,\dots) $ inductively by
\begin{align}
\notag
F_1(\gamma_1)&=\gamma_1, \quad F_2(\gamma_1,\gamma_2)=(1-|\gamma_1|^2)\gamma_2,\\
\label{eq:F}
F_{n}(\gamma_1,\dots, \gamma_{n})
&=(1-|\gamma_1|^2)F_{n-1}(\gamma_2,\dots, \gamma_{n}) \\
\notag
&\quad -\bar\gamma_1\sum_{k=2}^{n-1} F_{n-k}(\gamma_2,\dots,\gamma_{n-k+1})
F_k(\gamma_1,\dots, \gamma_k), \qquad n\ge 3.
\end{align}
By construction, the function $F_n(\gamma_1,\dots,\gamma_n)$
is a polynomial in $\gamma_1,\bar\gamma_1, \dots,$
$\gamma_{n-1},\bar\gamma_{n-1}, \gamma_n$
with integer coefficients.
From \eqref{eq:F} one easily infers that
\begin{align*}
F_3&=(1-|\gamma_1|^2)(1-|\gamma_2|^2)\gamma_3
-(1-|\gamma_1|^2)\bar\gamma_1\gamma_2^2, \\
F_4&=(1-|\gamma_1|^2)(1-|\gamma_2|^2)(1-|\gamma_3|^2)\gamma_4 \\
&\quad -(1-|\gamma_1|^2)(1-|\gamma_2|^2)\gamma_3(2\bar\gamma_1\gamma_2
+\bar\gamma_2\gamma_3)+(1-|\gamma_1|^2)\bar\gamma_1^2\gamma_2^3.
\end{align*}

By induction, the reader can easily verify the following Schur's recurrence relation.

\begin{lem}[\cite{Simon:OP1}]\label{lem:F}
For each $n\ge 2,$ there exists a function $G_n(\gamma_1,\dots,\gamma_{n-1})$
of $n-1$ complex variables $\gamma_1,\dots,\gamma_{n-1}$ such that the following
equality holds:
$$
G_n(\gamma_1,\dots,\gamma_{n-1})=F_n(\gamma_1,\dots,\gamma_n)-(1-|\gamma_1|^2)\cdots(1-|\gamma_{n-1}|^2)\gamma_n.
$$
\end{lem}

It is not difficult to prove the following result.

\begin{lem}[\cite{Li-Sugawa-Schur2019}]\label{lem:rec}
Let $g(z)=a_1z+a_2z^2+\cdots$ be a function in $\H_0$
with its hyperbolic derivatives $H^jg(0)=\gamma_j$ $(j=1,2,\ldots)$.
Then
$a_n=F_n(\gamma_1,\gamma_2,\dots, \gamma_n)$
for $n=1, 2, 3, \dots.$
\end{lem}

For $f\in \H$, Peschl's invariant derivatives $D^n f(z)$ are defined by the Taylor series expansion \cite{peschl1955invariants} (see also \cite{kim2007invariant} and \cite{Schippers07}):
$$z\mapsto g(z):=\frac{f(\dfrac{z+z_0}{1+\overline{z}_0 z})-f(z_0)}{1-\overline{f(z_0)}f(\dfrac{z+z_0}{1+\overline{z}_0 z})}=\sum_{n=1}^{\infty}\frac{D^n f(z_0)}{n!}z^n,\quad z, z_0\in \D,$$
where $D^n f(z_0)=g^{(n)}(0)$.
Precise forms of $D^n f(z)$, $n=1,2,3$, are expressed by
\begin{align*}
D_1 f(z)&=\frac{(1-\abs{z}^2)f'(z)}{1-\abs{f(z)}^2},\\
D_2 f(z)&=\frac{(1-\abs{z}^2)^2}{1-\abs{f(z)}^2}
\Bigg[f''(z)-\frac{2\overline{z}f'(z)}{1-\abs{z}^2}
+\frac{2\overline{f(z)}f'(z)^2} {1-\abs{f(z)}^2}\Bigg],\\
D_3 f(z)&=\frac{(1-\abs{z}^2)^3}{1-\abs{f(z)}^2}
\Bigg[f'''(z)-\frac{6\overline{z}f''(z)}{1-\abs{z}^2}
+\frac{6\overline{f(z)}f'(z)f''(z)} {1-\abs{f(z)}^2}
+\frac{6\overline{z}^2f'(z)}{(1-\abs{z}^2)^2}&\\
&\qquad\qquad\qquad -\frac{12\overline{f(z)}f'(z)^2} {(1-\abs{z}^2)(1-\abs{f(z)}^2)}
+\frac{6\overline{f(z)}^2f'(z)^3} {(1-\abs{f(z)}^2)^2}\Bigg].
\end{align*}
The first relations between $H^n f$ and $D^n f$ are given by
\begin{align*}
&H^{1} f(z)=D^{1} f(z),\\
&H^{2} f(z)=\frac{D^{2} f(z)}{2\left(1-\left|H^{1} f(z)\right|^{2}\right)},\\
&H^{3} f(z)=\frac{3 \overline{H^{1} f(z)} H^{2} f(z) D^{2} f(z)+D^{3} f(z)}{6-6\left|H^{1} f(z)\right|^{2}-3 \overline{H^{2} f(z)} D^{2} f(z)}.
\end{align*}
We recall that a Blaschke product of degree $n \in \N$  is a function of the form      $ B(z)=e^{i \theta}\prod\limits_{j=1}^{n}
        [z,z_j], \quad z, z_j\in \D, \theta \in \mathbb{R}.$
We denote the class of Blaschke products of degree $n$ by $\mathcal{B}_n$.

\begin{lem}[\cite{cho2012multi}]\label{prop:rel}
Let $f\in\H$ and $z_0\in\D.$
Define $g\in\H$
by $g(z)=[f([z,-z_0]),f(z_0)].$
If $H^ng(0)=
\gamma_n$ for $n=1,2,\dots,$
then  $H^nf(z_0)=\gamma_n$ for $n=1,2,\dots$.
\end{lem}
When we express $g$ by the series expansion
$g(z)=\sum_{n=1}^\infty a_n z^n,$ the first several $\gamma_j$'s
are given by
\begin{align*}
\gamma_1&= a_1, \\
\gamma_2&= \frac{a_2}{1-|a_1|^2}, \\
\gamma_3&= \frac{a_3(1-|a_1|^2)+\bar a_1a_2^2}{(1-|a_1|^2)^2-|a_2|^2}.
\end{align*}
Noting the relation $a_n=g^{(n)}(0)/n!=D_nf(z_0)/n!$, we can obtain the formula of $D_nf(z_0)$,
\begin{equation}\label{eq:D-f-recurrence}
D_nf(z_0)=n!(1-|\gamma_1|^2)\cdots(1-|\gamma_{n-1}|^2)\gamma_n
+n!G_n(\gamma_1,\dots,\gamma_{n-1}).
\end{equation}
We are able to obtain a formula  for  $H^nf$  in  terms  of  $D^nf$  and  its  lower-order hyperbolic derivatives  $H^1f,\ldots,H^{n-1}f$.

\begin{thm}
Let $n\geq 2$, $f\in \H\setminus \cup_{j=0}^{j=n-1}\mathcal{B}_j$ and $z\in \D$.
Then
\begin{equation}\label{eq:H-n-f-formula}
H^nf(z)=\frac{D_nf(z)-n!G_n(H^1f(z),\dots,H^{n-1}f(z))}
{n!(1-|H^1f(z)|^2)\cdots(1-|H^{n-1}f(z)|^2)}.
\end{equation}
\end{thm}
Using the formula given in \eqref{eq:H-n-f-formula}, it is easy to exhibit an explicit form of hyperbolic derivatives. For instance, $G_2(H^1f(z))=0$ and
$G_3(H^1f(z),H^2f(z))=-(1-|H^1f(z)|^2)\overline{H^1f(z)}H^2f(z)$, then we immediately obtain
$$ H^2f(z)=\frac{D_2f(z)-2!G_2(H^1f(z))}
{2!(1-|H^1f(z)|^2)} =\frac{D_2f(z)}{2(1-|H^1f(z)|^2)}$$
and
\begin{align*}
H^3f(z)&=\frac{D_3f(z)-3!G_3(H^1f(z),H^{2}f(z))}
{3!(1-|H^1f(z)|^2)(1-|H^{2}f(z)|^2)}\\
&=\frac{D_3f(z)+6(1-|H^1f(z)|^2)\overline{H^1f(z)}H^2f(z)}
{6(1-|H^1f(z)|^2)(1-|H^{2}f(z)|^2)}.
\end{align*}

It  is  possible  to  obtain  expressions  for  $D^nf$  in  terms  of  $f^{(n)}$  and  its  lower-order  derivatives,  as  well  as
derivatives  $D^1f,\ldots,D^nf$.

\begin{lem}[\cite{kim2007invariant},\text{ Corollary 7.5}]
Let $f:\mathbb{D}\to\mathbb{D}$ be  a  holomorphic  function.
Then
$$D^nf=\sum_{k=1}^n\alpha_{n,k}
\frac{(\overline{z})^{n-k}(1-|z|^2)^kf^{(k)}(z)}{1-|f(z)|^2}-\sum_{k=2}^nk!(-\overline{f(z)})^{k-1}A_{n,k}(D^1f,\ldots,D^{n-k+1}f),$$
where $A_{n, k}$ is the Bell polynomial given by
\begin{align*}
A_{n,k}(x_1,\ldots,x_{n-k+1})&:= \sum_{j\in I_k}{n! \over j_1!j_2!\cdots j_{n-k+1}!}
\left({x_1\over 1!}\right)^{j_1}\left({x_2\over 2!}\right)^{j_2}\cdots\left({x_{n-k+1} \over (n-k+1)!}\right)^{j_{n-k+1}} \\
&= n! \sum_{j\in I_k} \prod_{i=1}^{n-k+1} \frac{x_i^{j_i}}{(i!)^{j_i} j_i!},
\end{align*}
$I_k$ consists of all multi-indexes $(j_1,\ldots,j_{n-k+1})$ such that
\begin{align*}
j_1,\ldots,j_{n-k+1}\ge 0,\\
j_1+j_2+\cdots+j_{n-k+1}=k,\\
j_1+2j_2+\cdots+(n-k+1)j_{n-k+1}=n,
\end{align*}
$\alpha_{n,k}$ is defined by
$$\alpha_{n,k}:=
\left\{\begin{array}{cc}(-1)^{n-k}\dfrac{n!(n-1)!}{k!(k-1)!(n-k)!},&\quad if\quad 1 \leq k \leq n ; \\
0,&otherwise.\end{array}\right.$$
\end{lem}
The ordinary derivative $f^{(n)}$ can also be expressed  in terms of the invariant derivatives $D^kf (k = 1, 2, ... , n)$ as follows.
\begin{lem}[\cite{kim2007invariant},\text{ Corollary 7.6}]
Let $f:\mathbb{D}\to\mathbb{D}$ be  a  holomorphic  function.
Then
$$\frac{(1-|z|^2)^n}{1-|f(z)|^2}\frac{f^{(n)}(z)}{n!}
=\sum_{k=1}^n\binom{n-1}{k-1}(\bar{z})^{n-k}\cdot b_k,$$
where
$$b_k=\sum_{l=1}^k\frac{l!}{k!}(-\overline{f(z)})^{l-1}A_{k,l}(D^1f,\ldots,D^{k-l+1}f).$$
\end{lem}
Since $\alpha_{n,n}=1$, we have
\begin{align*}
D^nf(z)&=\frac{(1-|z|^2)^n f^{(n)}(z)}{1-|f(z)|^2}+s_{n-1}(z),
\end{align*}
where
$$s_{n-1}(z)= \sum_{k=1}^{n-1}\alpha_{n,k}\frac{(\overline{z})^{n-k}(1-|z|^2)^kf^{(k)}(z)}{1-|f(z)|^2}
-\sum_{k=2}^nk!(-\overline{f(z)})^{k-1}A_{n,k}(D^1f(z),\ldots,D^{n-k+1}f(z))$$
depends only on $D^1f(z),\ldots,D^{n-1}f(z)$ whenever $z$ and $f(z)$ are given. By the relation \eqref{eq:H-n-f-formula},
we note that $s_{n-1}(z)$ depends only on $H^1f(z),\ldots,H^{n-1}f(z)$ other than $z$ and $f(z)$. For instance,  $s_1(z)=-2H^1f(z)(\bar z-\overline{f(z)}H^1f(z))$.

Together with \eqref{eq:D-f-recurrence}, we have
$f^{(n)}(z_0)=c_n+\rho_n \gamma_n$,
where
\begin{equation}
\left\{
\begin{aligned}
c_n&=\frac{(1-|\gamma_0|^2)}{(1-|z_0|^2)^n}
[n!G_n(\gamma_1,\ldots,\gamma_{n-1})
-s_{n-1}(z_0)],\\
\rho_n&=\frac{n!\prod_{k=0}^{n-1}(1-|\gamma_k|^2)}{(1-|z_0|^2)^n}.
\end{aligned}
\right.
\end{equation}

For $\varepsilon\in \overline{\D}$,
let $f_{\gamma,\varepsilon} (z) =
T_{\gamma_0} (T_{-z_0}(z) T_{\gamma_1}(\cdots T_{-z_0}(z) T_{\gamma_{n-1}}(\varepsilon T_{-z_0}(z)) \cdots ))$.
In addition,
we set $G_1:=0$ and $s_0(z_0):=0$, then we obtain the following Schwarz-Pick Lemma of $n$-th derivative.
\begin{thm}[Schwarz-Pick Lemma of $n$-th derivative]
\label{thm:Schwarz-Pick-Lemma-of-n-th-derivative}
Let $n \in {\mathbb N}$, $z_0 \in {\mathbb D}$ and $\gamma  = (\gamma_0, \ldots , \gamma_{n-1} ) \in \overline{{\mathbb D}}^{n}$. Suppose that
$f \in {\mathcal H}$, $f(z_0)=\gamma_0$, $H^1f(z_0)=\gamma_1,
   \cdots, H^{n-1}f(z_0)=\gamma_{n-1}$.
 \begin{enumerate}
 \item
\label{boundary}
If $|\gamma_1|<1, \ldots , |\gamma_{j-1}|<1$, $|\gamma_j|=1$, $\gamma_{j+1}= \cdots = \gamma_{n-1}=0$
for some $j=1, \ldots , n-1$, then
$f^{(n)}(z_0)=c_n$ and
\begin{equation*}
f (z) =
  T_{\gamma_0} (T_{-z_0}(z) T_{\gamma_1}(\cdots T_{-z_0}(z) T_{\gamma_{j-1}}(\gamma_j T_{-z_0}(z))  \cdots )),
\end{equation*}
which is a Blaschke product of degree $j$.

\item
\label{interior}
If $|\gamma_1|<1, \ldots , |\gamma_{n-1}|<1$, then
 the region of values of $f^{(n)}(z_0)$ is the closed disk
$$
\overline{\D}(c_n, \rho_n)
=\{f^{(n)}_{\gamma,\varepsilon}(z_0):\varepsilon \in\overline{\D} \},
$$
and $f(z)$ is the form of
\begin{equation*}
f (z) =
T_{\gamma_0} (T_{-z_0}(z) T_{\gamma_1}(\cdots T_{-z_0}(z) T_{\gamma_{n-1}}(T_{-z_0}(z) f^* (z)) \cdots )),
\end{equation*}
where $f^* \in \mathcal{S}$.
Further, $f^{(n)}(z_0)\in \partial\D(c_n, \rho_n)$ if and only if
$f(z)=f_{\gamma,\varepsilon}(z)$ for some constant $\varepsilon \in \partial \D$.
\end{enumerate}

\end{thm}
It is worth pointing out that the Schwarz-Pick Lemma (Theorem A) is a simple corollary of Theorem \ref{thm:Schwarz-Pick-Lemma-of-n-th-derivative} for the case $n=1$.
The so-called ``
Schwarz-Pick Lemma of second derivative" (cf. \cite[Theorem 2.2]{chen2024second}), which describes the range of values of $f''(z_0)$ for $z_0\in\D$, can also be  directly derived from Theorem \ref{thm:Schwarz-Pick-Lemma-of-n-th-derivative} for the case $n=2$.
We denote
$c_2$ and $\rho_2$ by
\begin{equation}
\left\{
\begin{aligned}
c_2=c_2(z_0,\gamma_0,\gamma_1)&= \frac{2(1-|\gamma_0|^2)}{(1-|z_0|^2)^2}
(\overline{z_0}-\overline{\gamma_0}\gamma_{1})\gamma_{1};\\
\rho_2=\rho_2(z_0,\gamma_0,\gamma_1)&=\frac{2(1-|\gamma_0|^2)(1-|\gamma_1|^2)}
{(1-|z_0|^2)^2}.
\end{aligned}
\right.
\end{equation}

\begin{customthm}{C (Schwarz-Pick Lemma of Second Derivative)}\label{thm:g-second}
Let $z_0, \gamma_0\in \D$, and $ \gamma_1 \in \ov{\D}$.
Suppose that $f\in \H$,
$f(z_0)=\gamma_0$ and  $H^1 f(z_0)=\gamma_1$.
Set
\begin{align*}
 &f_{\gamma_1}(z) = T_{\gamma_0}\left(\gamma_1 T_{-z_0}(z)\right),\\
& f_{\gamma_1,\alpha}(z)=T_{\gamma_0}\left(T_{-z_0}(z) T_{\gamma_1}(\alpha T_{-z_0}(z))\right).
\end{align*}
\begin{enumerate}
\item If $|\gamma_1|=1$, then $f''(z_0)=c_2$ and
    $f(z)=f_{\gamma_1} (z)$.
\item If $|\gamma_1|<1$, then the region of values of $f''(z_0)$ is the closed disk
$$
\overline{\D}(c_2, \rho_2)
=\{f''_{\gamma_1,\alpha}(z_0):\alpha\in\overline{\D} \},
$$
and $f(z)$ is the form of
$T_{\gamma_0}\left(T_{-z_0}(z) T_{\gamma_1}(T_{-z_0}(z)f^*(z))\right)$,
where $f^*\in\S$.
Further, $g''(z_0)\in \partial\D(c_2, \rho_2)$ if and only if
$f(z)=f_{\gamma_1,\alpha}(z)$ for some constant $\alpha \in \partial \D$.
\end{enumerate}
\end{customthm}

\section{Dieudonn\'e's Lemma of $n$-th derivative}
 It is natural for us to further study the $n$-th order derivative $f^{(n)}$ of $f\in \H_0$. We can apply Theorem \ref{thm:Schwarz-Pick-Lemma-of-n-th-derivative} to formulate a Dieudonn\'e's Lemma of $n$-th derivative.
 In fact, the purpose of this section is to  deal with a variability region problem.
 \begin{problem}
Let $z_0, w_0\in \mathbb{D}$, $\gamma_1,..., \gamma_{n-1} \in \overline{\D}$ with $|w_0|=s<r=|z_0|$. Set $\gamma_0=w_0/z_0$. Denote $w_k\in\overline{\D}(c'_k,\rho'_k)$ by
$w_k=c'_k+\rho'_kz_0\gamma_k/r$ for $k=1,\ldots,n-1$, where
$$c'_1=\frac{w_0}{z_0}, \qquad \rho'_1=\frac{r^2-s^2}{r(1-s^2)},$$
and
\begin{equation}
c'_k=kc_{k-1}+k\rho_{k-1}\gamma_{k-1}+z_0c_k,\quad\rho'_k=r\rho_k,\quad \emph{for}~ k=2,\ldots,n-1.
\end{equation}
Suppose that $h\in\mathcal{H}_0$, $h(z_0) = w_0, h'(z_0)=w_1,...,h^{(n-1)}(z_0)=w_{n-1}$.
 Determine the region of values of $h^{(n)}(z_0)$.
 \end{problem}
 We just need to consider the case $\gamma_1,..., \gamma_{n-1} \in \D$.
 Assume that $f(z)=h(z)/z$, then $f\in \H$. It is easily seen that
$f(z_0)=\gamma_0$,  $H^1f(z_0)=\gamma_1,
   \cdots, H^{n-1}f(z_0)=\gamma_{n-1}$.
We let $h_{\gamma,\varepsilon}(z)=zf_{\gamma,\varepsilon}(z)$
and
denote  the class $\H_0 (z_0,w_0,\ldots,w_{n-1})$ by
$$\H_0 (z_0,w_0,\ldots,w_{n-1}) =  \{ h \in \H_0 : h(z_0) =w_0, h'(z_0)=w_1,\ldots, h^{(n-1)}(z_0)=w_{n-1}\}.$$
  From the proof of Theorem \ref{thm:Schwarz-Pick-Lemma-of-n-th-derivative} and the relations $h_{\gamma,\varepsilon}^{(k)}(z)=kf^{(k-1)}(z)+zf^{(k)}(z)$ for $k=1,\ldots,n$, it is easy to confirm that $h_{\gamma,\varepsilon}\in \mathcal{H}_0 (z_0,w_0,\ldots,w_{n-1})$ and obtain
\begin{align}
\label{eq:f-4}
h^{(n)}_{\gamma,\varepsilon}(z_0) =c'_n+\rho'_n\varepsilon,\nonumber
\end{align}
 where
 \begin{equation}
\left\{
\begin{aligned}
c'_n&=nc_{n-1}+n\rho_{n-1}\gamma_{n-1}+z_0c_n;\\
\rho'_n&=r\rho_n.
\end{aligned}
\right.
\end{equation}
 The closed disk $\overline{\D}(c'_n,\rho'_n)$ is covered since $\varepsilon \in \overline{\D}$ is arbitrary.

 We know that $h^{(n)}(z_0)\in \partial \D(c'_n, \rho'_n)$ if and only if $h(z)=zf(z)$, where $f$ is a Blaschke product of degree $n$ satisfying $f(z_0)=\gamma_0$,  $H^1f(z_0)=\gamma_1,
   \cdots, H^{n-1}f(z_0)=\gamma_{n-1}$. From the proof of Theorem \ref{thm:Schwarz-Pick-Lemma-of-n-th-derivative},
we can easily check that $f(z)=f_{\gamma,\varepsilon}(z)$, $\varepsilon\in \partial \D$.
For $\varepsilon\in \overline{\D}$,
let $h_{\gamma,\varepsilon} (z) =zf_{\gamma,\varepsilon} (z)=
T_{\gamma_0} (T_{-z_0}(z) T_{\gamma_1}(\cdots T_{-z_0}(z) T_{\gamma_{n-1}}(\varepsilon T_{-z_0}(z)) \cdots ))$.
\begin{thm}[Dieudonn\'e's Lemma of the $n$-th order]\label{thm:lemnth}
Let $z_0, w_0\in \mathbb{D}$, $\gamma_1,..., \gamma_{n-1} \in \ov{\D}$ with $|w_0|=s<r=|z_0|$.
Suppose that $h\in\mathcal{H}_0$, $h(z_0) = w_0, h'(z_0)=w_1,...,h^{(n-1)}(z_0)=w_{n-1}$.
Set $\gamma_0=w_0/z_0$.
\begin{enumerate}
\item If $|\gamma_1|<1, \ldots , |\gamma_{j-1}|<1$, $|\gamma_j|=1$, $\gamma_{j+1}= \cdots = \gamma_{n-1}=0$
for some $j=1, \ldots , n-1$, then
$h^{(n)}(z_0)=c'_n$ and
\begin{equation*}
h (z) =
  z T_{\gamma_0} (T_{-z_0}(z) T_{\gamma_1}(\cdots T_{-z_0}(z) T_{\gamma_{j-1}}(\gamma_j T_{-z_0}(z))  \cdots )).
\end{equation*}

\item If $|\gamma_1|<1$, $|\gamma_2|<1$,..., $|\gamma_{n-1}|<1$, then the region of values of $h^{(n)}(z_0)$ is the closed disk
$\overline{\D}(c'_n, \rho'_n)=\{h^{(n)}_{\gamma,\varepsilon}(z_0):\varepsilon \in\overline{\D}\}$
and $h(z)$ is the form of
\begin{equation*}
h (z) =
zT_{\gamma_0} (T_{-z_0}(z) T_{\gamma_1}(\cdots T_{-z_0}(z) T_{\gamma_{n-1}}(T_{-z_0}(z) h^* (z)) \cdots )),
\end{equation*}
where $h^* \in \mathcal{S}$.
Further, $h^{(n)}(z_0)\in \p\D(c'_n, \rho'_n)$ if and only if
$h(z)=zf_{\gamma,\varepsilon}(z)=z T_{\gamma_0}\left(T_{-z_0}(z) T_{\gamma_1}(T_{-z_0}(z) T_{\gamma_2}(\cdots T_{-z_0}(z)T_{\gamma_{n-1}}(\varepsilon T_{-z_0}(z))\cdots))\right)$ for some constant $\varepsilon \in \partial \D$.
\end{enumerate}
\end{thm}

It is worth noting that, Dieudonn\'e's Lemma (Theorem B) is a straightforward corollary of Theorem \ref{thm:lemnth} for the case $n=1$. Moreover, for the case $n=2$,
the above theorem yields a Dieudonn\'e's Lemma of the second order, which was earlist obtained by Rivard \cite{rivard2013application} (see also \cite{cho2012multi}).
Denote $c'_2$ and $\rho'_2$ by
\begin{equation*}
\left\{
\begin{aligned}
c'_2=c'_2(z_0,w_0,\gamma_1)&=\frac{2(r^2-s^2)}{r^2(1-r^2)^2}
\gamma_1(1-\frac{z_0\overline{w_0}}{\overline{z_0}}\gamma_1),\\ \rho'_2=\rho'_2(z_0,w_0,\gamma_1)&=\frac{2(r^2-s^2)}{r(1-r^2)^2}(1-|\gamma_1|^2).
\end{aligned}
\right.
\end{equation*}
\begin{customthm}{D (Dieudonn\'e's Lemma of Second Derivative)}
Let $z_0, w_0\in \mathbb{D}$, $\gamma_1\in \ov{\D}$ with $|w_0|=s<r=|z_0|$,
$$w_1= c'_1+\rho'_1\frac{r\gamma_1}{\overline{z_0}}
=\frac{w_0}{z_0}+\frac{r^2-s^2}{\overline{z_0}(1-s^2)}\gamma_1.$$
 Suppose that $h\in\mathcal{H}_0$, $h(z_0) = w_0$ and $h'(z_0)=w_1$.
Set $\gamma_0=w_0/z_0$, and
\begin{align*}
&h_{\gamma_1}(z)=z T_{\gamma_0}(\gamma_1 T_{-z_0}(z)),\\
&h_{\gamma_1,\alpha}(z)=z T_{\gamma_0}\big(T_{-z_0}(z) T_{\gamma_1}(\alpha T_{-z_0}(z))\big).
\end{align*}
\begin{enumerate}
\item If $|\gamma_1|=1$, then $h''(z_0)=c'_2$ and $h(z)=h_{\gamma_1}(z)$.
\item If $|\gamma_1|<1$, then the region of values of $h''(z_0)$ is the closed disk
$$\overline{\D}(c'_2, \rho'_2)
=\{h''_{\gamma_1,\alpha}(z_0):\alpha\in \overline{\D}\},
$$
and $h(z)$ is the form of
$$z T_{\gamma_0}\left(T_{-z_0}(z) T_{\gamma_1}(h^*(z) T_{-z_0}(z))\right),$$
where $h^*\in\S$.
Further, $h''(z_0)\in \partial\D(c'_2,\rho'_2)$ if and only if
$h(z)=h_{\gamma_1,\alpha}(z)$ for some constant $\alpha\in \partial \D$.
\end{enumerate}
\end{customthm}



\section*{Acknowledgement}
The author would like to express his deep gratitude to Prof. Toshiyuki Sugawa for his valuable comments and instructive suggestions. 

\bibliographystyle{srtnumbered}

\end{document}